\documentclass[10pt]{article}

\usepackage{latexsym,amsfonts}
\usepackage{amssymb,amsthm,upref,amscd}
\usepackage[T1]{fontenc}
\usepackage{times}

\newtheorem{Theorem}{Theorem}[section]
\newtheorem{definition}[Theorem]{Definition}
\newtheorem{lemma}[Theorem]{Lemma}
\newtheorem{proposition}[Theorem]{Proposition}
\newtheorem{corollary}[Theorem]{Corollary}
\newtheorem{remark}[Theorem]{Remark}
\newtheorem{Open Problem}[Theorem]{Open Problem}

\makeatletter \@addtoreset{equation}{section} \makeatother
\begin{document}

\title{\bf Existence of nontrivial solutions for   periodic   Schr\"odinger
equations with new   nonlinearities
 }

\author{Shaowei Chen \thanks{ E-mail address:
swchen6@163.com (Shaowei Chen)}    \quad \quad Dawei Zhang
  \\ \\
\small \small\it School of Mathematical Sciences, Huaqiao
University, \\
\small Quanzhou  362021,  China\\ }

\date{}
\maketitle
\begin{minipage}{13cm}
{\small {\bf Abstract:} We study the  Schr\"{o}dinger equation:
\begin{eqnarray}
- \Delta u+V(x)u+f(x,u)=0,\qquad u\in
H^{1}(\mathbb{R}^{N}),\nonumber
\end{eqnarray} where $V$ is periodic and
$f$ is periodic in the $x$-variables, $0$ is in a gap of the
spectrum of the operator $-\Delta+V$.   We prove that under some
new   assumptions for $f$,  this equation has a nontrivial
solution. Our assumptions for the nonlinearity $f$ are very weak
and greatly different
from the known assumptions  in the literature.   
}
\end{minipage}

\section{Introduction and statement of results}\label{diyizhang}
In this paper, we consider the following Schr\"{o}dinger equation:
\begin{equation}\label{e1}
- \Delta u+V(x)u+f(x,u)=0 ,\qquad u\in H^{1}(\mathbb{R}^{N}),
\end{equation}
where $N\geq1$. For $V$ and $f,$ we assume
\begin{description}
\item{$(\bf{v}).$}
 $V\in C( \mathbb{R}^{N}) $ is 1-periodic in $x_j$ for $j = 1,\cdots,N$,
  0 is in a spectral gap $(-\mu_{-1}, \mu_1)$ of  $-\Delta+V$
 and $-\mu_{-1}$ and $\mu_1$ lie in the essential spectrum of $-\Delta+V.$

Denote $$\mu_0 := \min\{\mu_{-1}, \mu_1\}.$$ \item{$(\bf{f_1}).$}
$f\in C( \mathbb{R}^{N}\times\mathbb{R}) $ is 1-periodic in $x_j$
for $j = 1,\cdots,N$.  And there exist constants $C>0$ and $2< p<
2^*$ such that
$$|f(x,t)|\leq C(1+|t|^{p-1}),\ \forall (x,t)\in\mathbb{R}^N\times\mathbb{R}$$
where $2^*:= \left\{
\begin{array}{l}
\frac{2N}{N-2}, \ N\geq 3\\
\infty,  \quad N=1,2.\\
\end{array} \right.$

\item{$(\bf{f_2}).$} The limit $\lim_{t\rightarrow 0}f( x,t)/t = 0
$ holds uniformly for $x\in\mathbb{R}^{N}$. And there there exists
$D>0$ such that
\begin{eqnarray}\label{nv99ufhf12}
\inf_{x\in\mathbb{R}^N,|t|\geq D}\frac{f(x,t)}{t}>
\max_{\mathbb{R}^N}V_-.
\end{eqnarray}
where $V_\pm(x)=\max\{\pm V(x),0\}$, $\forall x\in\mathbb{R}^N$.

 \item{$(\bf{f_3}).$} For any
$(x,t)\in\mathbb{R}^N\times\mathbb{R}$, $\widetilde{F}(x,t)\geq
0,$ where
\begin{eqnarray}
\widetilde{F}(x,t):=\frac{1}{2}tf(x,t)-F(x,t)
\nonumber\end{eqnarray} and $F(x,t)=\int^{t}_{0}f(x,s)ds$.

\item{$(\bf{f_4}).$} There exist $0<\kappa<D$ and
$\nu\in(0,\mu_0)$ such that, for every
$(x,t)\in\mathbb{R}^N\times\mathbb{R}$ with $|t|<\kappa$,
\begin{eqnarray}\label{bc77fdtdr}
|f(x,t)|\leq\nu |t| \end{eqnarray}  and for every
$(x,t)\in\mathbb{R}^N\times\mathbb{R}$ with $\kappa\leq|t|\leq D$,
\begin{eqnarray}\label{bv77fkkl} \widetilde{F}(x,t)> 0.
\end{eqnarray}
\end{description}

\begin{remark}\label{gc66tdrd}
 By the definitions of $F$ and $\widetilde{F}$, it is easy to
verify that, for all
$(x,t)\in\mathbb{R}^N\times(\mathbb{R}\setminus\{0\})$,
$$\frac{\partial}{\partial t}\Big(\frac{F(x,t)}{t^2}\Big)=\frac{2\widetilde{F}(x,t)}{t^3}.$$
Together with $f(x,t)=o(t)$ as $|t|\rightarrow 0$ and $(\bf f_3)$,
this implies that
\begin{eqnarray}\label{h77s5aa9}F(x,t)\geq 0\ \mbox{for all}\
(x,t)\in\mathbb{R}^N\times\mathbb{R}.\end{eqnarray}
\end{remark}

 A solution $u$ of (\ref{e1}) is called nontrivial if
$u\not\equiv0.$  Our main results are as follows:

\begin{Theorem}\label{th1}
Suppose $\bf (v)$,  and $\bf ( f_{1})-\bf ( f_{4}) $ are
satisfied. Then Eq.(\ref{e1}) has a nontrivial solution.
\end{Theorem}

Note that
\begin{description}
\item{$(\bf{f'_2}).$} The limits $\lim_{t\rightarrow 0}f( x,t)/t =
0 $ and $\lim_{|t|\rightarrow\infty}\frac{f( x,t)}{t} =+\infty$
hold uniformly for $x\in\mathbb{R}^{N}$.
\end{description}
implies $\bf (f_2)$. We have the following corollary

\begin{corollary}\label{vd5rdf}
Suppose $\bf (v)$,  $\bf ( f_{1})$, $\bf (f'_2)$, $\bf ( f_{3}) $,
and $\bf ( f_{4}) $  are satisfied. Then Eq.(\ref{e1}) has a
nontrivial solution.
\end{corollary}

It is easy to verify that the condition
\begin{description}
\item{$(\bf{f'_4}).$} $\widetilde{F}(x,t)> 0$ for every
$(x,t)\in\mathbb{R}^N\times\mathbb{R}$.
\end{description} and the assumption that $f(
x,t)/t \rightarrow 0 $ as $t\rightarrow0$ uniformly for
$x\in\mathbb{R}^{N}$   imply $\bf (f_3)$ and  $\bf (f_4)$.
Therefore, we have the following corollary:

\begin{corollary}\label{vd5rdf}
Suppose $\bf (v)$,  $\bf ( f_{1})$, $\bf (f_2)$,  and $\bf (
f'_{4}) $  are satisfied. Then Eq.(\ref{e1}) has a nontrivial
solution.
\end{corollary}

Semilinear Schr\"odinger equations with periodic coefficients have
attracted much attention in recent
 years due to its numerous  applications. One can see \cite{Ackermann}-\cite{Dingbook}, \cite{HKS},
 \cite{jean}-\cite{PMilan}, \cite{schechter}-\cite{YCD} and the references therein.  In
 \cite{alama}, the authors used the dual variational method to
 obtain  a nontrivial solution of (\ref{e1}) with
 $f(x,t)=\pm W(x)|t|^{p-2}t$ , where $W$ is a asymptotically periodic
 function. In \cite{TW}, Troestler and  Willem firstly obtained
 nontrivial solutions for (\ref{e1}) with $f$ is a $C^1$ function
 satisfying  the Ambrosetti-Rabinowitz condition:
 $$(AR)\qquad \mbox{there exists }\ \alpha>2\ \mbox{such that for every}\ u\neq 0,\ 0<\alpha G(x,u)\leq g(x,u)u,$$
 where $g(x,u)=-f(x,u)$ and $G(x,u)=-F(x,u),$  and
    $$\Big|\frac{\partial f(x,u)}{\partial u}\Big|\leq
 C(|u|^{p-2}+|u|^{q-2})$$ with  $2<p<q<2^*$. Then,
in \cite{KS}, Kryszewski and Szulkin developed
 some   infinite-dimensional linking theorems. Using these
 theorems, they improved  Troestler and  Willem's results and  obtained  nontrivial solutions for (\ref{e1}) with
 $f$ only satisfying $\bf (f_1)$ and the $(AR)$ condition. These generalized
linking theorems were also  used by Li and Szulkin to obtain
nontrivial solution for (\ref{e1}) under some asymptotically
linear assumptions for $f$ (see \cite{LiSzulkin}). In
\cite{pankov} (see also \cite{PMilan}), existence of nontrivial
solutions for (\ref{e1}) under $\bf (f_1)$ and the $(AR)$
condition was also obtained by  Pankov and  Pfl\"uger through
approximating  (\ref{e1}) by a sequence of equations defined in
bounded domains. In the celebrated paper  \cite{Zou}, Schechter
and Zou combined a generalized
  linking theorem with the monotonicity methods  of Jeanjean (see
  \cite{jean}). They obtained a nontrivial solution of (\ref{e1})
  when $f$ exhibts the critical growth. A similar  approach was applied by
  Szulkin and Zou to obtain homoclinic orbits of asymptotically
  linear Hamiltonian systems (see \cite{szulkin}). Moreover, in
  \cite{yanheng} (see also \cite{Dingbook}), Li and Ding obtained nontrivial solutions for
  (\ref{e1}) under
  some new superlinear assumptions on $f$ different from the classical
  $(AR)$ conditions.

Our assumptions on $f$ are very weak and greatly different from
the assumptions mentioned above.  In fact, our assumptions $\bf
(f_1)-\bf (f_4)$ do not involve the properties of $f$ at infinity.
It may be asymptotically linear growth at infinity, i.e.,
$\limsup_{|t|\rightarrow\infty}\frac{f(x,t)}{t}<+\infty$ or
superlinear growth at infinity as well, i.e.,
$\liminf_{|t|\rightarrow\infty}\frac{f(x,t)}{t}=+\infty$.

In this paper, we use the generalized linking theorem for a class
of parameter-dependent functionals  (see \cite[Theorem 2.1]{Zou}
or Proposition \ref{b99d443} in the present paper) to
 obtain a sequence of approximate solutions for (\ref{e1}). Then,
 we  prove that these
 approximate solutions are bounded in $L^\infty(\mathbb{R}^N)$ and $H^1(\mathbb{R}^N)$ (see Lemma \ref{bc77ftrg} and \ref{x4esdftrg}).
  Finally, using the concentration-compactness principle, we obtain
 a nontrivial solution of (\ref{e1}).

 \medskip

 \noindent{\bf Notation.} $B_r(a)$ denotes the  open ball of radius $r$ and center $a$.
For a Banach space $E,$ we denote the dual space of $E$ by $E'$,
and denote   strong and   weak convergence in $E$  by
$\rightarrow$ and $\rightharpoonup$, respectively. For $\varphi\in
C^1(E;\mathbb{R}),$ we denote the Fr\'echet derivative of
$\varphi$ at $u$ by $\varphi'(u)$.  The Gateaux derivative of
$\varphi$ is denoted by $\langle \varphi'(u), v\rangle,$ $\forall
u,v\in E.$ $L^p(\mathbb{R}^N)$ denotes the standard $L^p$ space
$(1\leq p\leq\infty)$, and $H^1(\mathbb{R}^N)$ denotes the
standard Sobolev space with norm
$||u||_{H^1}=(\int_{\mathbb{R}^N}(|\nabla u|^2+u^2)dx)^{1/2}.$  We
use $O(h)$, $o(h)$ to mean $|O(h)|\leq C|h|,$ $o(h)/|h|\rightarrow
0$ as $|h|\rightarrow 0$.

\section{Existence of approximate solutions for Eq.(\ref{e1})}\label{nvb6ftrf}

Under the  assumptions $\bf (v)$, $\bf (f_1)$, and $\bf (f_2)$,
the functional
\begin{equation}\label{mnncb66dtd}
\Phi(u)=\frac{1}{2}\int_{\mathbb{R}^{N}}\left\vert \nabla
u\right\vert ^{2}dx
+\frac{1}{2}\int_{\mathbb{R}^{N}}V(x)u^{2}dx+\int
_{\mathbb{R}^{N}}F(x,u)dx
\end{equation}
is of class $C^1$ on $X := H^1(\mathbb{R}^N )$, and the critical
points of $\Phi$ are weak solutions of (\ref{e1}).

 Assume that $\bf (v)$ holds, and let $S=-\Delta+V$ be the
self-adjoint operator acting on $L^2(\mathbb{R}^N)$ with domain
$D(S)=H^2(\mathbb{R}^N)$.  By virtue of $\bf (v)$,  we have the
orthogonal decomposition $$L^2=L^2(\mathbb{R}^N)=L^++L^-$$ such
that $S$ is negative (resp.positive) in $L^-$(resp.in $L^+$). As
in \cite[Section 2]{yanheng} (see also \cite[Chapter
6.2]{Dingbook}), let $X=D(|S|^{1/2})$ be equipped with the inner
product
$$(u,v)=(|S|^{1/2}u, |S|^{1/2}v)_{L^2}$$and norm $||u||=|||S|^{1/2}u||_{L^2}$,
 where $ (\cdot,\cdot)_{L^2}$ denotes the
inner product of $L^2$. From $\bf (v),$ $$X=H^1(\mathbb{R}^N)$$
with equivalent norms. Therefore, $X$  continuously embeds in $L^q
(\mathbb{R}^N)$ for all $2\leq q\leq2N/(N-2)$  if $N\geq 3$ and
for all $q\geq 2$  if $N=1,2$. In addition, we have the
decomposition
$$X=X^++X^-,$$ where $X^\pm=X\cap L^\pm$ is orthogonal with
respect to both $(\cdot,\cdot)_{L^2}$ and $(\cdot,\cdot)$.
Therefore, for every $u\in X$ , there is a unique decomposition
$$ u=u^++u^-,\ u^\pm\in X^\pm$$
 with $(u^+,u^-)=0$
and \begin{eqnarray}\label{buucytdrd} \int_{\mathbb{R}^N}|\nabla
u|^2dx+\int_{\mathbb{R}^N}V(x)u^2dx=||u^+||^2-||u^-||^2,\ u\in X.
\end{eqnarray} Moreover,
\begin{eqnarray}\label{bcvttdref}
&&\mu_{-1}||u^-||^2_{L^2}\leq||u^-||^2,\quad \forall u\in X,
\end{eqnarray} and\begin{eqnarray}\label{bcvttdref2}
&&\mu_{1}||u^+||^2_{L^2}\leq||u^+||^2,\quad \forall u\in X.
\end{eqnarray}
 The functional $\Phi$ defined by (\ref{mnncb66dtd}) can
be rewritten as
\begin{eqnarray}\label{vcrdfd}
 \Phi(u)=\frac{1}{2}(||u^+||^2-||u^-||^2)+\psi(u),
\end{eqnarray}
where $$\psi(u)=\int_{\mathbb{R}^N}F(x,u)  dx.$$

The above variational setting for the functional
(\ref{mnncb66dtd}) is standard. One can consult \cite{yanheng} or
\cite{Dingbook} for more details.

Let  $\{e^\pm_k\}$ be the   total orthonormal sequence in $X^\pm.$
Let $P:X\rightarrow X^-$, $Q:X\rightarrow X^+$ be the orthogonal
projections. We define
$$|||u|||=\max\Big\{||Pu||,\sum^{\infty}_{k=1}\frac{1}{2^{k+1}}|(Qu,e^+_k)|\Big\}$$
on $X.$
 The topology
generated by $|||\cdot|||$ is denoted by $\tau$, and all
topological notation related to it will include this symbol.

\begin{lemma}\label{jf777fyrrr}
Suppose that $\bf (v)$ holds. Then
\begin{description}
\item{$(\bf{a}).$} $\max_{\mathbb{R}^N} V_-\geq\mu_{-1},$
 where $\mu_{-1}$  is defined in $(\bf v)$. \item{$(\bf{b}).$} For any $C>\mu_{-1},$ there
exists $u_0\in X^-$ with $||u_0||=1$ such that
$C||u_0||_{L^2}>1.$
\end{description}
\end{lemma}
\noindent{\bf Proof.} $(\bf{a}).$ We  apply an indirect argument,
and assume by contradiction that $$c:=\max_{\mathbb{R}^N}
V_-<\mu_{-1}.$$ From assumption $\bf (v)$, $-\mu_{-1}$ is
 in the essential spectrum of the operator (with domain
$D(L)=H^2(\mathbb{R}^N)$) $$L = -\Delta+ V:
L^2(\mathbb{R}^N)\rightarrow L^2(\mathbb{R}^N)
$$ Then, by the Weyl's
criterion (see, for example, \cite[Theorem VII.12]{reed1} or
\cite[Theorem 7.2]{hislop}), there exists a sequence
$\{u_n\}\subset H^2(\mathbb{R}^N)$ with the properties that
$||u_n||_{L^2}=1$, $\forall n$ and $||-\Delta u_n+Vu_n+\mu_{-1}
u_n||_{L^2}\rightarrow 0$.

Since $\mu_{-1}>c=\max_{\mathbb{R}^N} V_-$, we deduce that
$-V_-(x)+\mu_{-1}> 0$ for all $x\in\mathbb{R}^N$. Together with
the facts that $V$ is a continuous periodic function and
$V=V_+-V_-$, this  implies
$$\inf_{x\in\mathbb{R}^N}(V(x)+\mu_{-1})>0.$$
It follows that there exists a constant $C'>0$ such that
\begin{eqnarray}\label{mcnn7yfyf}
\int_{\mathbb{R}^N}(|\nabla u|^2+(V(x)+\mu_{-1})u^2)dx\geq
C'||u||^2,\ \forall u\in X.
\end{eqnarray}
Note that
$$\int_{\mathbb{R}^N}(-\Delta u_n+V(x)u_n+\mu_{-1} u_n)u_ndx=\int_{\mathbb{R}^N}(|\nabla u_n|^2+(V(x)+\mu_{-1})u^2_n)dx.$$
Together with (\ref{mcnn7yfyf}) and the fact that $||-\Delta
u_n+Vu_n+\mu_{-1} u_n||_{L^2}\rightarrow 0$ and $||u_n||_{L^2}=1,$
this implies $||u_n||\rightarrow 0.$  It contradicts
$||u_n||_{L^2}=1,$ $\forall n.$ Therefore, $\max_{\mathbb{R}^N}
V_-\geq\mu_{-1}.$

$(\bf{b}).$ It suffices to prove that
$$\mu_{-1}=C_-:=\inf\{||u||^2\ |\ u\in X^-,\ ||u||_{L^2}=1\}.$$ From
(\ref{bcvttdref}), we deduce that $\mu_{-1}\leq C_-.$ From
assumption $\bf (v)$, $-\mu_{-1}$ is  in the essential spectrum of
$L$. By the Weyl's criterion, there exists $\{u_n\}\subset
H^2(\mathbb{R}^N)$ such that $||u_n||_{L^2}=1$ and $||-\Delta
u_n+Vu_n+\mu_{-1} u_n||_{L^2}\rightarrow 0$. Multiplying  $-\Delta
u_n+Vu_n+\mu_{-1} u_n$ by $u^+_n$ and then integrating on
$\mathbb{R}^N$, by (\ref{buucytdrd}) and (\ref{bcvttdref2}), we
get that
\begin{eqnarray}(\mu_1+\mu_{-1})||u^+_n||^2_{L^2}&\leq&
\int_{\mathbb{R}^N}(|\nabla
u^+_n|^2+V(x)(u^+_n)^2+\mu_{-1}(u^+_n)^2)dx\nonumber\\
&=&\int_{\mathbb{R}^N}(-\Delta u_n+V(x)u_n+\mu_{-1}
u_n)u^+_ndx\rightarrow 0.\nonumber \end{eqnarray} It follows that
$||u^-_n||_{L^2}\rightarrow1$. Multiplying $-\Delta
u_n+Vu_n+\mu_{-1} u_n$ by $u^-_n$ and
 then integrating on $\mathbb{R}^N$, we get that
\begin{eqnarray}-||u^-_n||^2+\mu_{-1}||u^-_n||^2_{L^2}
&=&\int_{\mathbb{R}^N}(|\nabla
u^-_n|^2+V(x)(u^-_n)^2+\mu_{-1}(u^-_n)^2)dx\nonumber\\
&=&\int_{\mathbb{R}^N}(-\Delta u_n+Vu_n+\mu_{-1}
u_n)u^-_ndx\rightarrow 0. \end{eqnarray} It implies that
$\mu_{-1}\geq C_-.$ This together with $\mu_{-1}\leq C_-$  implies
$\mu_{-1}= C_-$. \hfill$\Box$

\medskip

Let $R > r > 0$ and $$A:=\inf_{x\in\mathbb{R}^N,|t|\geq
D}\frac{f(x,t)}{t}.$$ From  assumption (\ref{nv99ufhf12}), we have
$A>\max_{\mathbb{R}^N}V_-$. Together with the result $\bf (a)$ of
Lemma \ref{jf777fyrrr}, this implies that
$\frac{1}{2}(A+\mu_{-1})>\mu_{-1}$. Choose
\begin{eqnarray}\label{n88cuyddd}
\gamma\in (\mu_{-1},(A+\mu_{-1})/2).
\end{eqnarray}
Then by the result $\bf (b)$ of Lemma \ref{jf777fyrrr},  there
exists $u_0 \in X^-$ with $||u_0||
 = 1$  such that \begin{eqnarray}\label{nvbyyftfg}
\gamma ||u_0||_{L^2}>1.
 \end{eqnarray}
 Set
$$N = \{u \in X^- \ |\  ||u||
 = r\},\  M = \{u \in X^+ \oplus \mathbb{R}^+u_0 \ |\
||u|| \leq R\}.$$ Then, $M$ is a submanifold of $X^+ \oplus
\mathbb{R}^+u_0$ with boundary \begin{eqnarray}
\partial M=\{u\in X^-\ |\ ||u||\leq
R\}\cup\{u^-+tu_0\ |\ u^-\in X^-,\ t>0,\ ||u^-+tu_0||=R
\}.\nonumber
\end{eqnarray}

\begin{definition}\label{bcv99iuytg}
Let  $\phi\in  C^1(X; \mathbb{R})$. A sequence $\{u_n\} \subset X$
is called a Palais-Smale sequence at  level $c$ ($(PS)_c$-sequence
for short) for $\phi$, if $\phi(u_n) \rightarrow c$ and
$||\phi'(u_n)||_{X'}\rightarrow 0$ as $n\rightarrow\infty.$
\end{definition}

The following proposition is proved in  \cite{Zou} (see
\cite[Theorem 2.1]{Zou}).
\begin{proposition}
 \label{b99d443} Let $0<K<1.$ The family
of $C^1$-functional $\{H_\lambda\}$ has the form
\begin{eqnarray}\label{bcttfrd}
H_\lambda(u)=\lambda I(u)- J(u), \ u\in X,\ \lambda\in[K,1].
\end{eqnarray}
Assume
\begin{description} \item {$(a)$} $J(u)\geq 0$, $\forall u\in X$, \item
{$(b)$} $|I(u)|+J(u)\rightarrow+\infty$ as
$||u||\rightarrow+\infty$, \item {$(c)$} for all
$\lambda\in[K,1]$, $H_\lambda$ is $\tau$-sequentially upper
semi-continuous, $i.e.,$ if $|||u_n-u|||\rightarrow 0,$ then
$$\limsup_{n\rightarrow\infty}H_\lambda(u_n)\leq H_\lambda(u),$$
and $H'_\lambda$  is weakly sequentially continuous. Moreover, $
H_\lambda$ maps bounded sets to bounded sets, \item {$(d)$} there
exist $u_0 \in X^- \setminus\{0\} $ with $||u_0||=1$, and  $R > r
> 0$ such that for all $\lambda\in[K,1]$,
$$\inf_N H_\lambda>\sup_{\partial M} H_\lambda.$$
\end{description} Then there exists $E\subset[K,1]$ such that the Lebesgue measure
of $[K,1]\setminus E$ is zero and for every  $\lambda\in E$, there
exist $c_\lambda$  and a bounded $(PS)_{c_\lambda}$-sequence
 for $H_\lambda,$ where $c_\lambda$ satisfies
\begin{eqnarray}\label{bvzcx}
\sup_{M}H_\lambda\geq c_\lambda\geq\inf_N
H_\lambda.\nonumber\end{eqnarray}
\end{proposition}

\medskip

For $0<K<1$ and  $\lambda\in[K,1]$, let
\begin{eqnarray}\label{vc6drsd}
 \Psi_\lambda(u)=\frac{\lambda}{2}\int_{\mathbb{R}^N}V_-(x)u^2dx-\Big(\frac{1}{2}\int_{\mathbb{R}^N}(|\nabla
 u|^2+V_+(x)u^2)dx+\psi(u)\Big),\
 u\in X.
\end{eqnarray}
Then $\Psi_1=-\Phi$ and it is easy to verify that a critical point
$u$ of $\Psi_\lambda$ is a weak solution of
\begin{eqnarray}\label{bv66frd}
-\Delta u+V_\lambda(x)u +f(x,u)=0,\ u\in X,
\end{eqnarray}
where $$V_\lambda=V^+-\lambda V^-.$$

\begin{lemma} \label{nxc5xrdf}
Suppose that $\bf (v)$ and $\bf (f_1)-(f_3)$  hold. Then, there
exist $0<K_*<1$ and  $E\subset[K_*,1]$ such that the Lebesgue
measure of $[K_*,1]\setminus E$ is zero and, for every $\lambda\in
E$, there exist $c_\lambda$ and a bounded
$(PS)_{c_\lambda}$-sequence for $\Psi_\lambda,$ where $c_\lambda$
satisfies
$$+\infty>\sup_{\lambda\in E}c_\lambda\geq\inf_{\lambda\in
E}c_\lambda>0.$$
\end{lemma}
\noindent{\bf Proof.} For $u\in X,$ let
$$I(u)=\frac{1}{2}\int_{\mathbb{R}^N}V_-(x)u^2dx$$
and $$J(u)=\frac{1}{2}\int_{\mathbb{R}^N}(|\nabla
u|^2+V_+(x)u^2)dx+\psi(u).$$ Then,  $I$ and $J$ satisfy
assumptions $(a)$ and $(b)$ in Proposition \ref{b99d443}, and, by
(\ref{vc6drsd}), $\Psi_\lambda(u)=\lambda I(u)- J(u)$.

From (\ref{vc6drsd}) and (\ref{buucytdrd}),   for any $u\in X$ and
$\lambda\in[K,1]$, we have
\begin{eqnarray}\label{bc77dtdf}
\Psi_{\lambda}(u)&=&\frac{\lambda-1}{2}\int_{\mathbb{R}^{N}}V_-(x)u^2dx-\Big(\frac{1}{2}\int_{\mathbb{R}^{N}}(\left\vert
\nabla u\right\vert ^{2}+V(x)u^{2})dx+\int
_{\mathbb{R}^{N}}F(x,u)dx\Big)\nonumber\\
&=&\frac{1}{2}||u^-||^2-\frac{1}{2}||u^+||^2-\frac{1-\lambda}{2}\int_{\mathbb{R}^{N}}V_-(x)u^2dx-\int
_{\mathbb{R}^{N}}F(x,u)dx.
\end{eqnarray}
Let $u_*\in X$ and $\{u_n\}\subset X$ be such that
$|||u_n-u_*|||\rightarrow 0$. It follows that $u_n^-\rightarrow
u^-_*$, $u^+_n\rightharpoonup u^+_*$, and $u_n\rightharpoonup
u_*$. In addition, up to a subsequence, we can assume that
$u_n\rightarrow u_*$ $a.e.$ in $\mathbb{R}^N$. Then, we have
\begin{eqnarray}
&&||u^-_n||^2\rightarrow ||u^-_*||^2,\nonumber\\
&&\liminf_{n\rightarrow\infty}\int_{\mathbb{R}^{N}}V_-(x)u^2_ndx\geq
\int_{\mathbb{R}^{N}}V_-(x)u^2_*dx\quad (\mbox{by the Fatou's
lemma}
),\nonumber\\
&&\liminf_{n\rightarrow\infty}||u^+_n||^2\geq
||u^+_*||^2.\nonumber\
\end{eqnarray}
By Remark \ref{gc66tdrd},  $F(x,t)\geq 0$ for all $x$ and $t.$
This together with the Fatou's lemma implies
\begin{eqnarray}
\liminf_{n\rightarrow\infty}\int _{\mathbb{R}^{N}}F(x,u_n)dx\geq
\int _{\mathbb{R}^{N}}F(x,u_*)dx.\nonumber
\end{eqnarray}
Then, by  (\ref{bc77dtdf}), we obtain
\begin{eqnarray}\label{nnvb77dtdf}
\limsup_{n\rightarrow\infty}\Psi_\lambda(u_n)\leq\Psi_\lambda(u_*).\nonumber
\end{eqnarray}
This implies that $\Psi_\lambda$ is $\tau$-sequentially upper
semi-continuous.

If $u_n\rightharpoonup u_*$ in $X,$ then, for any fixed
$\varphi\in X$, as $n\rightarrow\infty,$
\begin{eqnarray}\label{mmbn88dfd}
\langle-\Psi'_\lambda(u_n),\varphi\rangle
&=&\int_{\mathbb{R}^N}(\nabla u_n\nabla \varphi+V_\lambda
u_n\varphi)dx+\int_{\mathbb{R}^N}f(x,u_n)\varphi
dx\nonumber\\
&\rightarrow&\int_{\mathbb{R}^N}(\nabla u_*\nabla
\varphi+V_\lambda u_*\varphi)dx+\int_{\mathbb{R}^N}f(x,u_*)\varphi
dx\nonumber\\
&=&\langle-\Psi'_\lambda(u_*),\varphi\rangle.\nonumber
\end{eqnarray}
This implies that $\Psi'_\lambda$ is weakly sequentially
continuous. Moreover, it is easy to see that $\Psi_\lambda$ maps
bounded sets to bounded sets. Therefore, $\Psi_\lambda$ satisfies
 assumption $(c)$  in Proposition \ref{b99d443}.

Finally, we shall verify    assumption $(d)$  in Proposition
\ref{b99d443} for $\Psi_\lambda$.

From assumption $\bf (f_1)$ and  $f(x,t)/t\rightarrow 0$ as
$t\rightarrow 0$ uniformly for $x\in\mathbb{R}^N$, we deduce that
for any $\epsilon>0,$ there exists $C_\epsilon>0$ such that
\begin{eqnarray}\label{m8bugygyyy}
F(x,t)\leq \epsilon t^2+C_\epsilon|t|^p,\ \forall (x,t)\in
\mathbb{R}^N\times \mathbb{R}.
\end{eqnarray}
From (\ref{bc77dtdf}) and (\ref{m8bugygyyy}),  we have, for $u\in
N,$
\begin{eqnarray}
\Psi_\lambda(u)\geq\frac{1}{2}||u||^2-\frac{1-\lambda}{2}\int_{\mathbb{R}^{N}}V_-(x)u^2dx-\epsilon\int
_{\mathbb{R}^{N}}u^2dx-C_\epsilon\int
_{\mathbb{R}^{N}}|u|^pdx.\nonumber
\end{eqnarray}
Then by the Sobolev inequality $||u||_{L^p(\mathbb{R}^N)}\leq
C||u||$ and $||u||_{L^2}\leq C||u||$ (by  (\ref{bcvttdref}) and
(\ref{bcvttdref2})), we deduce that there exists a constant $C>0$
such that
\begin{eqnarray}
\Psi_\lambda(u)\geq\frac{1}{2}||u||^2-C(1-\lambda)\max_{\mathbb{R}^N}V_-(x)||u||^2-\epsilon
C||u||^2-CC_\epsilon||u||^p.\nonumber
\end{eqnarray}
Choose $0<K_*<1$ and $\epsilon>0$ such that
$C(1-K_*)\max_{\mathbb{R}^N}V_-(x)<1/4$ and $C\epsilon=1/8$. Then
for every $\lambda\in [K_*,1]$,  we have
\begin{eqnarray}\label{12mmbn8ftfg}
\Psi_\lambda(u)\geq\frac{1}{8}||u||^2-CC_{\epsilon}||u||^p.
\end{eqnarray}
Let  $r>0$ be such that $r^{p-2}CC_\epsilon=1/16$ and
$\beta=r^2/16.$  Then from (\ref{12mmbn8ftfg}), we deduce that,
for $ N=\{u\in X^-\ |\ ||u||=r\},$
\begin{eqnarray}\label{oiifyuft66}
\inf_{N}\Psi_\lambda \geq\beta,\ \forall\lambda\in
[K_*,1].\nonumber
\end{eqnarray}

We shall prove that $\sup_{K_*\leq\lambda\leq
1}\Psi_{\lambda}(u)\rightarrow-\infty$ as $||u||\rightarrow\infty$
and  $u\in X^+\oplus \mathbb{R}^+u_0$. Arguing indirectly, assume
that for some sequences $\lambda_n\in[K_*,1]$ and $u_n\in
X^+\oplus \mathbb{R}^+u_0$ with $||u_n||\rightarrow+\infty$, there
is $ \mathcal{L}
> 0$ such that $\Psi_{\lambda_n}(u_n)\geq-\mathcal{L}$ for all $n$. Then, setting
$w_n=u_n/||u_n||$,   we have $||w_n||=1$, and, up to a
subsequence, $w_n\rightharpoonup w$, $w^-_n\rightarrow w^-\in X^-$
and $w^+_n\rightharpoonup w^+\in X^+.$

First, we consider the case $w\neq 0.$ Dividing both sides of
(\ref{bc77dtdf}) by $||u_n||^2$, we get that
\begin{eqnarray}\label{mmvn77dtdg}
-\frac{\mathcal{L}}{||u_n||^2}\leq\frac{\Psi_{\lambda_n}(u_n)}{||u_n||^2}
=\frac{1}{2}||w_n^-||^2-\frac{1}{2}||w^+_n||^2-\frac{1-\lambda_n}{2}\int_{\mathbb{R}^{N}}V_-(x)w_n^2dx-\int
_{\mathbb{R}^{N}}\frac{F(x,u_n)}{||u_n||^2}dx.
\end{eqnarray}

From  (\ref{nv99ufhf12}) and the result $\bf (a)$ of Lemma
\ref{jf777fyrrr}, we deduce that
\begin{eqnarray}\label{nv99viuuvj}
\liminf_{|t|\rightarrow\infty}\frac{F(x,t)}{t^2}\geq
\frac{A}{2}>\frac{1}{2}\max_{\mathbb{R}^N}V_-\geq\frac{1}{2}\mu_{-1},\nonumber
\end{eqnarray}
where $A:=\inf_{x\in\mathbb{R}^N,|t|\geq D}\frac{f(x,t)}{t}$. Note
that for $x\in\left\{  x\in\mathbb{R}^{N}\mid w\neq0\right\} $, we
have $| u_{n}(x)| \rightarrow+\infty$. This implies that, when $n$
is large enough,
$$\int_{\{
x\in\mathbb{R}^{N}\mid w\neq0
\}}\frac{F(x,u_n)}{u_n^2}w^2_ndx\geq\frac{A+\mu_{-1}}{4}\int_{\{
x\in\mathbb{R}^{N}\mid w\neq0 \}}w^2_ndx.$$  By (\ref{h77s5aa9}),
we have, when $n$ is large enough,
\begin{eqnarray} \int
_{\mathbb{R}^{N}}\frac{F(x,u_n)}{||u_n||^2}dx=\int
_{\mathbb{R}^{N}}\frac{F(x,u_n)}{u_n^2}w^2_ndx\geq \int _{\{
x\in\mathbb{R}^{N}\mid w\neq0
\}}\frac{F(x,u_n)}{u_n^2}w^2_ndx.\nonumber
\end{eqnarray}  Combining the above two inequalities yields
\begin{eqnarray}\label{nnvb77dyd}
&&\liminf_{n\rightarrow\infty}\Big(\frac{1}{2}||w_n^-||^2-\frac{1}{2}||w^+_n||^2
-\frac{1-\lambda_n}{2}\int_{\mathbb{R}^{N}}V_-(x)w_n^2dx-\int
_{\mathbb{R}^{N}}\frac{F(x,u_n)}{||u_n||^2}dx\Big)
\nonumber\\
&\leq&\liminf_{n\rightarrow\infty}\Big(\frac{1}{2}||w_n^-||^2-\frac{1}{2}||w^+_n||^2-\frac{A+\mu_{-1}}{4}\int_{\{
x\in\mathbb{R}^{N}\mid w\neq0 \}}w^2_ndx\Big)\nonumber\\
&\leq&\frac{1}{2}||w^-||^2-\frac{1}{2}||w^+||^2-\frac{A+\mu_{-1}}{4}\int_{\mathbb{R}^{N}}w^2dx\nonumber\\
&\leq&\frac{1}{2}||w^-||^2-\frac{1}{2}||w^+||^2-\frac{A+\mu_{-1}}{4}||w^-||^2_{L^2}.
\end{eqnarray}
We used the inequalities
$$\lim_{n\rightarrow\infty}||w^-_n||^2= ||w^-||^2,\ \liminf_{n\rightarrow\infty}||w^+_n||^2\geq ||w^+||^2\ \mbox{ and}\
\liminf_{n\rightarrow\infty}\int_{\{ x\in\mathbb{R}^{N}\mid w\neq0
\}}w^2_ndx\geq \int_{\mathbb{R}^{N}}w^2dx$$ in the second
inequality of (\ref{nnvb77dyd}).

 Since $w^-=tu_0$ for some $t\in\mathbb{R}$, by (\ref{nvbyyftfg}), we get that
$$\frac{A+\mu_{-1}}{4}||w^-||^2_{L^2}\geq \frac{A+\mu_{-1}}{4\gamma}||w^-||^2.$$
Note that, by the choice of $\gamma$ (see (\ref{n88cuyddd})), we
have $\frac{A+\mu_{-1}}{4\gamma}>1/2$. Then  by (\ref{nnvb77dyd})
and the fact that $w\neq0$, we have that
\begin{eqnarray}\label{jfgfggf998dss}
&&\liminf_{n\rightarrow\infty}\Big(\frac{1}{2}||w_n^-||^2-\frac{1}{2}||w^+_n||^2
-\frac{1-\lambda_n}{2}\int_{\mathbb{R}^{N}}V_-(x)w_n^2dx-\int
_{\mathbb{R}^{N}}\frac{F(x,u_n)}{||u_n||^2}dx\Big)
\nonumber\\
&\leq&-\Big(\frac{A+\mu_{-1}}{4\gamma}-\frac{1}{2}\Big)||w^-||^2-\frac{1}{2}||w^+||^2<0.\nonumber
\end{eqnarray}
It contradicts (\ref{mmvn77dtdg}), since
$-\mathcal{L}/||u_n||^2\rightarrow0$ as $n\rightarrow\infty$.

Second, we consider the case $w= 0.$ In this case,
$\lim_{n\rightarrow\infty}||w^-_n||=0.$ It follows that
$$\liminf_{n\rightarrow\infty}||w^+_n||\geq1,$$ since $||w_n||=1$ and
$w_n=w^+_n+w^-_n$. Therefore, the right hand side of
(\ref{mmvn77dtdg}) is less than $-1/4$ when $n$ is large enough.
However, as $n\rightarrow\infty$, the left hand side of
(\ref{mmvn77dtdg}) converges to zero. It   induces a
contradiction.

Therefore, there exists $R>r$ such that
$$\sup_{\lambda\in[K_*,1]}\sup_{\partial M}\Psi_\lambda\leq 0.$$
This implies that $\Psi_\lambda$ satisfies  assumption $(d)$ in
Proposition \ref{b99d443} if $\lambda\in[K_*,1]$. Finally, it is
easy to see that
$$\sup_{\lambda\in[K_*,1]}\sup_{M}\Psi_\lambda<+\infty.$$
 Then, the results of this lemma follow immediately
from Proposition \ref{b99d443}. \hfill$\Box$

\medskip

\begin{lemma}\label{bc6dr5df} Suppose that $\bf (v)$ and $\bf (f_1)-(f_3)$  are
satisfied.   Let $\lambda\in[K_*,1]$ be fixed, where $K_*$ is  the
constant  in Lemma \ref{nxc5xrdf}. If $\{v_{n}\}$ is a bounded
$(PS)_c$ sequence for $\Psi_\lambda$ with $c\neq0$, then for every
$n\in\mathbb{N},$ there exists $a_n\in\mathbb{Z}^N$ such that, up
to a subsequence, $u_n:=v_n(\cdot+a_n)$ satisfies
\begin{eqnarray}
u_n\rightharpoonup u_\lambda\neq 0,\quad
 \Psi_\lambda(u_\lambda)\leq c \quad \mbox{and}\quad
\Psi'_\lambda(u_\lambda)= 0.
\end{eqnarray}
\end{lemma}
\noindent{\bf Proof.} The proof of this lemma is inspired by the
proof of Lemma 3.7 in \cite{szulkin}. Because $\{v_n\}$ is a
bounded sequence in $X,$   up to a subsequence, either
\begin{description}
\item {$(a)$}
$\lim_{n\rightarrow\infty}\sup_{y\in\mathbb{R}^N}\int_{B_1(y)}|v_n|^2dx=0$,
or \item {$(b)$} there exist $\varrho>0$ and $a_n\in\mathbb{Z}^N$
such that $\int_{B_1(a_n)}|v_n|^2dx\geq\varrho.$
\end{description}

If  $(a)$ occurs,  using the Lions lemma (see, for example,
\cite[Lemma 1.21]{Willem}), a similar argument as for the proof of
\cite[Lemma 3.6]{szulkin} shows that
\begin{eqnarray}\label{vvcgtdrdf}
\lim_{n\rightarrow\infty}\int_{\mathbb{R}^N}F(x,v_n)dx=0\quad
\mbox{and}\quad
\lim_{n\rightarrow\infty}\int_{\mathbb{R}^N}f(x,v_n)v^\pm_ndx=0.
\end{eqnarray}
It follows that \begin{eqnarray}\label{bcv66ftheshu}
\lim_{n\rightarrow\infty}\int_{\mathbb{R}^N}(2F(x,v_n)-f(x,v_n)v_n)dx=0.
\end{eqnarray} On the other hand, as $\{v_n\}$ is a $(PS)_c$ sequence of
$\Psi_\lambda,$ we have $\langle\Psi'_\lambda(v_n),
v_n\rangle\rightarrow 0$ and $\Psi_\lambda(v_n)\rightarrow c\neq
0$. It follows that
\begin{eqnarray}
&&\int_{\mathbb{R}^N}(f(x,v_n)v_n-2F(x,v_n))dx\nonumber\\
&=&2\Psi_\lambda(v_n)-\langle\Psi'_\lambda(v_n),
v_n\rangle\rightarrow 2c\neq0,\quad n\rightarrow\infty.
\end{eqnarray}
This contradicts (\ref{bcv66ftheshu}). Therefore,  case $(a)$
cannot occur.

If case  $(b)$ occurs, let $u_n=v_n(\cdot+a_n)$.   For every $n,$
\begin{eqnarray}\label{ggcfdrdf}
\int_{B_1(0)}|u_n|^2dx\geq\varrho.
\end{eqnarray}
Because $V$ and $F(x,t)$ are $1$-periodic in every $x_j$,
$\{u_n\}$ is still bounded in $X$,
\begin{eqnarray}\label{bbc66ftr} \lim_{n\rightarrow\infty}\Psi_\lambda(u_n)\leq
c\quad \mbox{and}\quad \Psi'_\lambda(u_n)\rightharpoonup0,\quad
n\rightarrow\infty.
\end{eqnarray}   Up to a subsequence, we  assume that
$u_n\rightharpoonup u_\lambda$ in $X$ as $n\rightarrow\infty$.
Since $u_n \rightarrow u_\lambda$ in $L^2_{ loc}( \mathbb{R}^N)$,
it follows from (\ref{ggcfdrdf}) that $u_\lambda \neq0$. Recall
that $\Psi'_\lambda(u_n)$ is weakly sequentially continuous.
Therefore,
$\Psi'_\lambda(u_n)\rightharpoonup\Psi'_\lambda(u_\lambda)$ and,
by (\ref{bbc66ftr}), $\Psi'_\lambda (u_\lambda)=0.$

Finally, by $\bf(f_3)$ and the Fatou's lemma
\begin{eqnarray}
c&=&\lim_{n\rightarrow\infty}(\Psi_\lambda(u_n)-\frac{1}{2}\langle\Psi'_\lambda(u_n),u_n\rangle)\nonumber\\
&=&\lim_{n\rightarrow\infty}\int_{\mathbb{R}^N}\widetilde{F}(x,u_n)
\geq\int_{\mathbb{R}^N}\widetilde{F}(x,u_\lambda)=\Psi_\lambda(u_\lambda).\nonumber
\end{eqnarray}\hfill$\Box$

\medskip

\begin{lemma}\label{ncbvdrdf}
There exist $0<K_{**}<1$ and $\eta>0$ such that for any
$\lambda\in[K_{**},1]$, if $u\neq 0$ satisfies
$\Psi'_\lambda(u)=0,$ then $||u||\geq\eta.$
\end{lemma}
\noindent{\bf Proof.} We adapt the arguments of Yang \cite[ p.
2626]{yang} and Liu \cite[Lemma 2.2]{liu}.
   Note that by $\bf (f_1)$ and $\bf (f_2)$,
for any $\epsilon > 0$, there exists $C_\epsilon > 0$ such that
$$|f(x,t)|\leq \epsilon |t|+C_\epsilon|t|^{p-1}.$$
Let $u\neq 0$ be a critical point of $\Psi_\lambda$. Then  $u$ is
a solution of
$$-\Delta u+V_{\lambda}u+f(x,u)=0,\ u\in X.$$
Multiplying both sides of this equation by $u^\pm$ respectively
and then integrating on $\mathbb{R}^N$, we get that
\begin{eqnarray}\label{ncb66ftf}
0&=&\pm||u^\pm||^2+(1-\lambda)\int_{\mathbb{R}^N}V_-(x)u_nu^\pm
dx+\int_{\mathbb{R}^N}f(x,u )u^\pm dx.\nonumber
\end{eqnarray}
It follows that
\begin{eqnarray}\label{kvxcdxssv}
||u^\pm||^2&=&\mp(1-\lambda)\int_{\mathbb{R}^N}V_-(x)uu^\pm
dx\mp\int_{\mathbb{R}^N}f(x,u)u^\pm dx\\
&\leq&(1-\lambda)\sup_{\mathbb{R}^N}V_-\int_{\mathbb{R}^N}|u|\cdot|u^\pm|
dx\nonumber\\
&&+\epsilon\int_{\mathbb{R}^N}|u|\cdot|u^\pm|
dx+C_\epsilon\int_{\mathbb{R}^N}|u|^{p-1}|u^\pm|dx\nonumber\\
&\leq&C_1((1-\lambda)+\epsilon)||u||\cdot||u^\pm||+C_2||u||^{p-1}||u^\pm||,\nonumber
\end{eqnarray}
where $C_1$ and $C_2$ are positive constants related to the
Sobolev inequalities, and $\sup_{\mathbb{R}^N}V_-.$ From the above
two inequalities, we obtain
\begin{eqnarray}
||u||^2=||u^+||^2+||u^-||^2\leq
2C_1((1-\lambda)+\epsilon)||u||^2+2C_2||u||^{p}.\nonumber
\end{eqnarray}
Because $p > 2$, this implies that $||u||
 \geq\eta$ for some $\eta> 0$ if  $\epsilon>0$ and $1-K_{**}>0$ are small enough
 and $\lambda\in[K_{**},1]$. The desired result follows.
 \hfill$\Box$

\medskip

 Let $K=\max\{K_*,K_{**}\}$, where
$K_*$ and $K_{**}$ are the constants  that appeared  in Lemma
\ref{nxc5xrdf} and Lemma \ref{ncbvdrdf}, respectively. Combining
Lemmas $\ref{nxc5xrdf} - \ref{ncbvdrdf}$, we obtain the following
lemma:
\begin{lemma}\label{bv88f7f} Suppose $\bf (v)$  and $\bf ( f_{1})-\bf ( f_{3}) $  are satisfied. Then, there exist $\eta>0,$
 $\{\lambda_n\}\subset[K,1]$, and $\{u_n\}\subset X$ such
that $\lambda_n\rightarrow 1$,
$$\sup_{n}\Psi_{\lambda_n}(u_n)<+\infty, \quad ||u_n||\geq\eta, \quad
\mbox{and}\quad \Psi'_{\lambda_n}(u_n)=0.$$
\end{lemma}

\section{A priori bound of approximate solutions  and proof of the main Theorem}

In this section, we give a priori bound for   the sequence of
approximate solutions $\{u_n\}$ obtained in Lemma \ref{bv88f7f}.
We then give the proofs of Theorem \ref{th1}.

\begin{lemma}\label{bc77ftrg} Suppose $\bf (v)$,  and $\bf ( f_{1})-\bf ( f_{3}) $  are satisfied.
 Let $\{u_n\}$ be the sequence
obtained in Lemma \ref{bv88f7f}. Then, $\{u_n\}\subset
L^\infty(\mathbb{R}^N)$ and
\begin{eqnarray}\label{bc77cyfdt}
\sup_{n}||u_n||_{L^\infty(\mathbb{R}^N)}\leq D.
\end{eqnarray}
\end{lemma}
\noindent{\bf Proof.} From $\Psi'_{\lambda_n}(u_n)=0$, we deduce
that $u_n$ is a weak solution of (\ref{bv66frd}) with
$\lambda=\lambda_n,$ i.e.,
\begin{eqnarray}\label{bv77ftfg}
-\Delta u_n+V_{\lambda_n}(x)u_n+f(x,u_n)=0\quad\mbox{in}\quad
\mathbb{R}^N.
\end{eqnarray}
By assumption $\bf (f_1)$ and  the bootstrap argument of elliptic
equations, we  deduce  that $u_n\in L^\infty(\mathbb{R}^N)$.

Multiplying both sides of (\ref{bv77ftfg}) by $v_n=(u_n-D)^+$ and
integrating on $\mathbb{R}^N$, we get that
\begin{eqnarray}\label{nnviiokjj}
\int_{\mathbb{R}^N}|\nabla v_n|^2dx+ \int_{u_n\geq
D}(V_{\lambda_n}(x)u_n+f(x,u_n))v_ndx=0.
\end{eqnarray}
Recall  that $V_{\lambda_n}=V^+-\lambda_nV^-$ and $\lambda_n\leq
1$. Then by (\ref{nv99ufhf12}), we get that
\begin{eqnarray}
 \int_{u_n\geq
D}(V_{\lambda_n}(x)u_n+f(x,u_n))v_ndx= \int_{u_n\geq
D}\Big(V_{\lambda_n}(x)+\frac{f(x,u_n)}{u_n}\Big)u_nv_ndx\geq
0.\nonumber
\end{eqnarray}
This together with (\ref{nnviiokjj}) yields $v_n=0.$ It follows
that $u_n(x)\leq D$  on $\mathbb{R}^N$.

Similarly,  multiplying both sides of (\ref{bv77ftfg}) by
$w_n=(u_n+D)^-$ and integrating on $\mathbb{R}^N$, we can get that
$u_n\geq-D$  on $\mathbb{R}^N$. Therefore, for all $n,$
$||u_n||_{L^\infty(\mathbb{R}^N)}\leq D$.\hfill$\Box$

\begin{lemma}\label{x4esdftrg} Suppose that $\bf (v)$, $\bf ( f_{1})$, $\bf (f_2)$, $\bf ( f_{3}) $, and
$\bf ( f_{4}) $  are satisfied. Let $\{u_n\}$ be the sequence
obtained in Lemma \ref{bv88f7f}. Then
\begin{eqnarray}\label{bc777eedt}
0<\inf_{n}||u_n||\leq\sup_{n}||u_n||<+\infty.
\end{eqnarray}
\end{lemma}
\noindent{\bf Proof.} As $\Psi'_{\lambda_n}(u_n)=0$ and $u_n\neq
0,$   Lemma \ref{ncbvdrdf} implies that $\inf_{n}||u_n||>0.$

To prove $\sup_{n}||u_n||<+\infty$, we apply an indirect argument,
and assume by contradiction that $||u_n||\rightarrow+\infty.$

Since $\Psi'_{\lambda_n}(u_n)=0$,  by (\ref{kvxcdxssv}), we get
that
\begin{eqnarray}\label{ncb66ftf}
||u^\pm_n||^2&=&\mp(1-\lambda_n)\int_{\mathbb{R}^N}V_-(x)u_nu^\pm_ndx\mp\int_{\mathbb{R}^N}f(x,u_n)u^\pm_ndx\nonumber\\
&=&\mp\int_{\mathbb{R}^N}f(x,u_n)u^\pm_ndx+(1-\lambda_n)O(||u_n||^2).\nonumber
\end{eqnarray}
It follows that
\begin{eqnarray}\label{hf55drzm}
&&||u_n||^2+\int_{\mathbb{R}^N}f(x,u_n)(u^+_n-u^-_n)dx\nonumber\\
&=&||u^+_n||^2+||u^-_n||^2+\int_{\mathbb{R}^N}f(x,u_n)(u^+_n-u^-_n)dx=(1-\lambda_n)O(||u_n||^2).
\end{eqnarray}
Set $w_n=u_n/||u_n||$. Then by (\ref{hf55drzm}),
$$||u_n||^2\Big(1+\int_{\mathbb{R}^N}\frac{f(x,u_n)}{u_n}(w^+_n-w^-_n)w_ndx\Big)=(1-\lambda_n)O(||u_n||^2).$$
Then by $\lambda_n\rightarrow 1$ as $n\rightarrow\infty,$ we have
that
\begin{eqnarray}\label{bv88dtdf}
\int_{\mathbb{R}^N}\frac{f(x,u_n)}{u_n}(w^+_n-w^-_n)w_ndx\rightarrow
-1,\quad n\rightarrow\infty.
\end{eqnarray}

From Lemma \ref{bv88f7f},
$$C_0:=\sup_{n}\Psi_{\lambda_n}(u_{n})<+\infty.$$
Then, by $\Psi'_{\lambda_n}(u_n)=0$, we obtain
\begin{eqnarray}\label{bc66drdf}
2C_0&\geq&2\Psi_{\lambda_n}(u_{n})-\langle
\Psi'_{\lambda_n}(u_{n}),
u_n\rangle=2\int_{\mathbb{R}^N}\widetilde{F}(x,u_n)dx\nonumber
\end{eqnarray}
From  $\bf (f_3)$, we have
\begin{eqnarray}\label{nb8fdr5er}
2C_0\geq 2\int_{\mathbb{R}^N}\widetilde{F}(x,u_n)dx\geq
2\int_{\{x\ |\ D\geq|u_n(x)|\geq\kappa\}}\widetilde{F}(x,u_n)dx
\end{eqnarray}
where $\kappa$ is the constant in $\bf (f_4)$.  As the continuous
function $\widetilde{F}$ is $1$-periodic in every $x_j$ variable,
we deduce from (\ref{bv77fkkl}) that there exists a constant
$C'>0$ such that
\begin{eqnarray}\label{bcv66fgff}
\widetilde{F}(x,t)\geq C't^2,\quad \mbox{for every }\
(x,t)\in\mathbb{R}^N\times\mathbb{R} \  \mbox{with}\
\kappa\leq|t|\leq D,
\end{eqnarray}
Combining (\ref{nb8fdr5er}) and (\ref{bcv66fgff}) leads to
$$C_0\geq C'\int_{\{x\ |\
D\geq|u_n(x)|\geq\kappa\}}u^2_ndx.$$ Dividing both sides of this
inequality by $||u_n||^2$ and sending $n\rightarrow\infty$, we
obtain
\begin{eqnarray}\label{bbvc77ftgd}
\lim_{n\rightarrow\infty}\int_{\{x\ |\
D\geq|u_n(x)|\geq\kappa\}}w^2_ndx=0.
\end{eqnarray}

From (\ref{bc77fdtdr}), (\ref{bcvttdref}), and (\ref{bcvttdref2}),
we have that
\begin{eqnarray}\label{nncb77dtdf}
&&\int_{\{x\ |\
|u_n(x)|<\kappa\}}\Big|\frac{f(x,u_n)}{u_n}(w^+_n-w^-_n)w_n\Big|dx\nonumber\\
&\leq&\nu\int_{\{x\ |\
|u_n(x)|<\kappa\}}|(w^+_n-w^-_n)w_n|dx\nonumber\\
&\leq&\nu\int_{\mathbb{R}^N}|(w^+_n-w^-_n)w_n|dx\nonumber\\
&\leq&\nu||w_n||^2_{L^2}\leq\frac{\nu}{\mu_0}||w_n||^2=\frac{\nu}{\mu_0}<1
\end{eqnarray}
where $\mu_0$ is the constant defined in $\bf (v)$.

Since $f\in C(\mathbb{R}^N\times \mathbb{R})$ and
$\lim_{t\rightarrow 0}f(x,t)/t=0,$ we deduce that there exists
$C>0$ such that for every $(x,t)\in\mathbb{R}^N\times\mathbb{R}$
with $|t|\leq D,$
$$|f(x,t)|\leq C|t|.$$ This together with (\ref{bbvc77ftgd}) gives
\begin{eqnarray}\label{nv88ihj}
&&\int_{\{x\ |\
D\geq|u_n(x)|\geq\kappa\}}\Big|\frac{f(x,u_n)}{u_n}(w^+_n-w^-_n)w_n\Big|dx\nonumber\\
&\leq&C\int_{\{x\ |\
D\geq|u_n(x)|\geq\kappa\}}|(w^+_n-w^-_n)w_n|dx\nonumber\\
&\leq& C||w^+_n-w^-_n||_{L^2}\Big(\int_{\{x\ |\
D\geq|u_n(x)|\geq\kappa\}}w^2_ndx\Big)^{1/2}\nonumber\\
&\leq&C||w_n||_{L^2}\Big(\int_{\{x\ |\
D\geq|u_n(x)|\geq\kappa\}}w^2_ndx\Big)^{1/2}\rightarrow 0,\quad
n\rightarrow\infty.
\end{eqnarray}
Combining (\ref{nncb77dtdf}) and (\ref{nv88ihj}) yields
\begin{eqnarray}\label{nnv8yydfd}
&&\limsup_{n\rightarrow\infty}\int_{\mathbb{R}^N}\Big|\frac{f(x,u_n)}{u_n}(w^+_n-w^-_n)w_n\Big|dx\nonumber\\
&\leq&\limsup_{n\rightarrow\infty}\int_{\{x\ |\
|u_n(x)|<\kappa\}}\Big|\frac{f(x,u_n)}{u_n}(w^+_n-w^-_n)w_n\Big|dx\nonumber\\
&&+\limsup_{n\rightarrow\infty}\int_{\{x\ |\
D\geq|u_n(x)|\geq\kappa\}}\Big|\frac{f(x,u_n)}{u_n}(w^+_n-w^-_n)w_n\Big|dx<1.\nonumber
\end{eqnarray}
This contradicts (\ref{bv88dtdf}). Therefore, $\{u_n\}$ is bounded
in $X.$\hfill$\Box$

\bigskip

\noindent{\bf Proof of Theorem \ref{th1}.} Let $\{u_n\}$ be the
sequence obtained in Lemma \ref{bv88f7f}.  From Lemma
\ref{x4esdftrg},  $\{u_n\}$ is bounded in $X$. Therefore,  up to a
subsequence, either
\begin{description}
\item {$(a)$}
$\lim_{n\rightarrow\infty}\sup_{y\in\mathbb{R}^N}\int_{B_1(y)}|u_n|^2dx=0$,
or \item {$(b)$} there exist $\varrho>0$ and $y_n\in\mathbb{Z}^N$
such that $\int_{B_1(y_n)}|u_n|^2dx\geq\varrho.$
\end{description}
According to (\ref{vvcgtdrdf}), if  case $(a)$ occurs,
\begin{eqnarray}
\lim_{n\rightarrow\infty}\int_{\mathbb{R}^N}f(x,u_n)u^\pm_ndx=0.\nonumber
\end{eqnarray}
Then, by (\ref{kvxcdxssv}) and $\lambda_n\rightarrow 1$, we have
\begin{eqnarray}\label{gdf55ere}
||u^\pm_n||^2&=&\mp(1-\lambda_n)\int_{\mathbb{R}^N}V_-(x)u_nu^\pm_n
dx\mp\int_{\mathbb{R}^N}f(x,u_n)u^\pm_n dx\nonumber\\
&\leq&
C(1-\lambda_n)||u_n||^2_{L^2}+\Big|\int_{\mathbb{R}^N}f(x,u_n)u^\pm_n
dx\Big|\rightarrow 0.
\end{eqnarray}
This contradicts $\inf_{n}||u_n||>0$ (see (\ref{bc777eedt})).
Therefore,  case $(a)$ cannot occur. As case  $(b)$ therefore
occurs, $w_n=u_n(\cdot+y_n)$ satisfies $w_n\rightharpoonup u_0\neq
0$. From  (\ref{mnncb66dtd}) and (\ref{vc6drsd}), we have that
$$\Psi_{\lambda}(u)=-\Phi(u)+\frac{\lambda-1}{2}\int_{\mathbb{R}^N}V_-u^2dx,\ \forall u\in X.$$
It follows that
\begin{eqnarray}\label{mc99c9c97}
\langle\Psi'_{\lambda}(u),\varphi\rangle=-\langle\Phi'(u),\varphi\rangle+(\lambda-1)\int_{\mathbb{R}^N}V_-u\varphi
dx,\ \forall u,\varphi\in X.
\end{eqnarray}By  $\Psi'_{\lambda_n}(u_n)=0$ (by Lemma
\ref{bv88f7f}), we have $\Psi'_{\lambda_n}(w_n)=0$. From
(\ref{mc99c9c97}), we have that, for any $\varphi\in X,$
\begin{eqnarray}\label{bv88ftdfd}
\langle\Psi'_{\lambda_n}(w_n),\varphi\rangle=-\langle\Phi'(w_n),\varphi\rangle+(\lambda_n-1)\int_{\mathbb{R}^N}V_-(x)w_n\varphi
dx.\nonumber
\end{eqnarray}
Together with  $\Psi'_{\lambda_n}(w_n)=0$ and
$\lambda_n\rightarrow 1$, this yields
$$\langle\Phi'(w_n),\varphi\rangle\rightarrow 0,\quad \forall\varphi\in X.$$
Finally, by $w_n\rightharpoonup u_0\neq 0$ and the weakly
sequential continuity of $\Phi'$, we have that $\Phi'(u_0)=0.$
Therefore, $u_0$ is a nontrivial solution of Eq.(\ref{e1}). This
completes the proof.\hfill$\Box$

\bigskip

{\bf \Large Conflict of Interests}

 \medskip

 The authors declare that there is no
conflict of interests regarding the publication of this paper.

\end{document}